\documentclass[3p, 11pt]{elsarticle}
\usepackage{lineno}
\usepackage[english]{babel}
\usepackage{natbib}
\usepackage{amssymb}
\usepackage{setspace}
\usepackage{lipsum, color}
\usepackage{amsmath}
\usepackage{natbib}
\usepackage{multirow}
\usepackage{longtable}
\usepackage{float}
\usepackage{adjustbox}
\usepackage{cuted}
\usepackage{float}
\usepackage{chngcntr}
\usepackage{mathptmx}
\usepackage{setspace}
\usepackage{isomath}
\usepackage{subfig}
\usepackage{textcomp}
\usepackage{bm}
\usepackage{lineno}

\numberwithin{equation}{section}

\newtheorem{propos}{Proposition}[section]

\biboptions{authoryear}

\makeatletter
\def\ps@pprintTitle{%
   \let\@oddhead\@empty
   \let\@evenhead\@empty
   \let\@oddfoot\@empty
   \let\@evenfoot\@oddfoot
}
\makeatother

\begin{document}

\begin{frontmatter}

\title{An outranking Choquet integral formulation for multi-criteria sorting problems with heterogeneous scales: an extension of FlowSort for interacting criteria}


\author[rvt]{Renata~Pelissari\corref{cor1}}
\ead{renatapelissari@unicamp.br}

\author[rvt]{Leonardo~Tomazeli~Duarte}
\ead{leonardo.duarte@fca.unicamp.br}

\cortext[cor1]{Corresponding author}

\address[rvt]{School of Applied Sciences, UNICAMP, Pedro Zaccaria 1300, 13484-350 Limeira, SP, Brazil}

\begin{abstract}
In multi-criteria decision aiding, the Choquet integral has been used as an aggregation operator to deal with interacting criteria, having as a requirement a prior assumption of common scale for all the criteria. This restriction on the adopted scale may be considered as a limitation in some practical problems. In order to overcome this limitation, we propose a new Choquet integral formulation, specific for sorting problems, that constructs a common scale from heterogenous scales using the framework of the FlowSort method. We also introduce a new outranking degree based on the Choquet integral preference model, which allows the modeling of interacting criteria in FlowSort. Therefore, the proposed approach can be seen either as an extension of FlowSort for problems with interacting criteria or as a new Choquet integral formulation for multi-criteria sorting problems with heterogeneous scales. Numerical examples attest that the proposed approach is conceptually simple to be implemented and acknowledge the importance of considering interaction between criteria when it exists. 
\end{abstract}


\end{frontmatter}

\onehalfspacing

\section{Introduction}

    \noindent In multi-criteria decision aiding (or multi-criteria decision making--MCDA/M),  there are four different kinds of analyses that can be performed in order to provide significant support to decision-makers \cite{book.roy1996}: (1) to identify the major distinguishing features of the alternatives and describe them based on these features (descriptive problem), (2) to identify the best alternatives (choice problem), (3) to construct a ranking of the alternatives from best to worst (ranking problem), (4) to sort the alternatives into predefined and ordered categories (sorting problem). Our interest in this paper is in problems of type 4.
    
    In mathematical terms, the problem addressed here concerns the sorting of $m$ alternatives $A=\{a_1, \ldots, a_m\}$ to $k$ predefined ordered categories $K_1, K_2, \ldots, K_k$. The alternatives are evaluated with respect to a set of criteria $G = \{g_1, g_2, \ldots, g_n\}$, and the performance evaluations of the alternatives, denoted by $g_j(a_i)$, $j=1, \ldots, n$ and $i=1, \ldots, m$, are the elements of the decision matrix presented in \eqref{decision_matrix}, whose rows are given by the alternatives of set $A$ and columns are the criteria of $G$:

    \begin{equation}\label{decision_matrix}
        M = \left[
            \begin{array}{cccc}
            g_1(a_1) & g_2(a_1) & \ldots & g_n(a_1) \\
            g_1(a_2) & g_2(a_2) & \ldots & g_n(a_2) \\
            \vdots & \vdots & \ldots & \vdots       \\
            g_1(a_m) & g_2(a_m) & \ldots & g_n(a_m)
            \end{array}
            \right].
    \end{equation}
    The categories may be defined through limiting profiles, which indicate the performance levels on each criterion that separate the categories. 
     
    We consider here that the criteria within $G$ may have different importance weights and there may also be interaction between any two of them. For instance, let us consider a decision problem of buying a vehicle, in which set $A$ is composed of different car models while criteria in $G$ are features of the cars such as price, maximum speed, acceleration, comfort, and so on. There may be an interaction between the criteria maximum speed and acceleration, since, in general, speedy cars also have good acceleration. Thus, their comprehensive importance is smaller than the sum of the importance of the two criteria considered separately. This type of interaction is called redundancy. 
    
    Analogously, the criteria comfort and price lead to a synergy effect, another type of interaction, since a comfortable car also having a low price is very well appreciated. Thus, the comprehensive importance of these two criteria should be greater than the sum of their importance considered separately. 
    
    In order to represent preferences in case of interaction between criteria, one has to use more general preference models than the weighted sum. This is the case of the non-additive integrals, among which the most well-known example is the Choquet integral \citep{Choquet_1954}. The Choquet integral is an extension of the weighted sum model to the case of interacting criteria, and it is based on the concept of capacity (fuzzy measure) that assigns a weight to each subset of criteria \citep{GRABISCH1996445, Grabisch2008}.
    
    Generally, the Choquet integral is applied directly to the decision matrix $M$ and a fundamental requirement is that the performance evaluation of the alternatives must be on the same scale for all criteria. Indeed, the Choquet integral preference model is based on the comparison between the performance evaluation of alternatives, and therefore the use of different scales is not possible. This requirement is known as the commensurability assumption \citep{Modave1998}.
    
    A problem that arises in the application of the Choquet integral is that the commensurability assumption is not an easy requirement to address. This requirement is easily met when the criteria can naturally be evaluated with the same quantitative scale, as in the case of evaluation of students regarding their grades in different subjects, where the same scale of grades can naturally be assumed for all subjects. When the criteria are evaluated by quantitative heterogeneous scales, a way to meet the commensurability assumption is by adopting a qualitative scale or degrees of preference to represent the heterogeneous scales \citep{Modave1998}. Adopting this strategy, we end up with a new limitation where one may consider that assuming qualitative scales, when quantitative information is available, might result in a loss of information. 
    
    Another way to meet the commensurability assumption, when heterogeneous scales are used, is by applying normalization techniques in order to convert the heterogeneous scales into a common scale, usually between 0 and 1. Two of the most applied normalization procedures in different fields (not only with the Choquet integral preference model) are min-max and z-score, which use maximum and minimum values, and mean and standard deviation concepts, respectively, in the scale conversion \citep{book-han2006}.
    
    In the MCDA/M context, \cite{Angilella2015a} proposed a normalization procedure specifically for the Choquet integral preference model based on a uniform sampling in the range [0,1], in which the generated values replace the heterogeneous scale. A drawback of this procedure is that it implies we shall assume the data are uniformly distributed, which is not always true. This procedure may also be considered inappropriate since the difference between two values is important in the Choquet integral, and it is lost when a uniform sample is used as a replacement for the heterogeneous scales.
    
    Our aim here is to propose a new Choquet integral formulation, in the context of MCDA/M, for sorting problems without any prior assumption about commensurability. The proposed method is based on the construction of common scales by applying the concept of preference function and flow used in the methods of the PROMETHEE family (see \cite{VinckeBrans1985, BEHZADIAN2010198} to know more about PROMETHEE). Particularly, in our proposal, we use the framework of FlowSort, an extension of PROMETHEE for sorting problems \citep{Nemery200890}.
    
    In the proposed method, pairwise comparisons (between the evaluation of the alternatives and the limiting profiles that define the categories) are made using preference functions which lead us to a common scale. From that moment on, the evaluations defined by heterogeneous scales are no longer used. Then, the flows are computed using the Choquet integral instead of the weighted sum. Sorting rules similar to the ones of FlowSort are applied to sort the alternatives into categories. It is worth noticing that the proposed method follows a similar idea to that presented by  \cite{Corrente2014}, in which the bi-Choquet integral \citep{GRABISCH2005211, GRABISCH2005237, GRECO201321} is integrated with PROMETHEE. 
    
    An advantage of the proposed formulation is that neither prior assumption about commensurability nor data uniformly distributed is required. Another benefit of the proposed method is to take advantage of the FlowSort framework, since it allows the use of preference functions to take the DM preference into account. Finally, as the proposed method is based on FlowSort, it is a non-compensatory method.
    
    It is worth emphasizing that the proposed method can also be seen as an extension of the FlowSort method for interacting criteria. Although there are some extensions of FlowSort to deal with interval data \citep{Janssen2013171}, qualitative data \citep{campos2015}, uncertainty data and preference elicitation \citep{PELISSARI2019235}, there is no version of FlowSort to model interacting criteria. Thereby, the contributions of the method proposed here can be summarized in two points: (i) the possibility of applying the Choquet integral for sorting problems using heterogeneous scales on the criteria evaluation; and (ii) an extension of FlowSort for interacting criteria. 

    This paper is organized as follows. In Section 2, we introduce some basic concepts required for the development of the proposed method, starting with the introduction of the Choquet integral and the problem of dealing with heterogeneous scales. Then we briefly describe the FlowSort method. In Section 3, we introduce the proposed method, which in this paper is designed FlowSort-Choquet, and we also present some of its proprieties in Section 4. In Section 5, the proposed method is illustrated by a numerical example, and in Section 6 we conduct some numerical tests. Conclusions and future research trends are presented in Section 7.
    
\section{Preliminaries}
    \noindent In this section, we introduce some basic concepts that are the background for the method proposed in this paper. 
    
    \subsection{The Choquet integral preference model} \label{secao_choquet_integral}
        
            \noindent Let be a decision-making problem with a set $A=\{a_1, \ldots, a_m\}$ of $m$ alternatives evaluated with respect to a set of $n$ criteria $G=\{g_1, g_2, \ldots, g_n\}$. The Choquet integral is based on the concept of capacity (fuzzy measure) that assigns a weight to each subset of criteria. More precisely, denoting by $2^G$ the power set of $G$, the function $\mu:2^G \rightarrow [0,1]$ is called a capacity on $2^G$ if the following conditions are satisfied:
            \begin{itemize}
                \item [(a)] Boundary conditions: $\mu(\emptyset) = 0$ and $\mu(G) = 1$,
                \item [(b)] Monotonicity condition: $\forall S \subseteq T \subseteq G, \mu(S) \leq \mu(T)$.
            \end{itemize}
            
            In a decision problem where the criteria are interacting and their importance is represented by a capacity $\mu$, the weighted sum can be extended through the Choquet integral which is a function $CI: A \rightarrow \mathbb{R}^{+}_{0}$ given by
            \begin{equation} \label{eq_choquet}
                CI(a) = \sum_{j=1}^{n}[g_{(j)}(a) - g_{(j-1)}(a)] \mu (N_j),
            \end{equation}
            where $g_{(0)}(a) \leq g_{(1)}(a) \leq \ldots \leq g_{(n)}(a)$, $g_{(0)}(a)=0$, $N_j=\{(j), \ldots, (n)\}$, for $j=1, \ldots, n.$
            
            A meaningful and useful reformulation of the capacity $\mu$ can be obtained by the M\"obius representation which is a function $m:2^{G} \rightarrow \mathbb{R}$ \citep{Shafer-1976} defined by
            \begin{equation}
                \mu(S) = \sum_{T \in S} m(T).
            \end{equation}
            
            In terms of M\"obius representation, boundary and monotonicity conditions presented in (a) and (b) are reformulated to the conditions presented in (a.2) and (b.2), respectively
             \begin{itemize}
                \item [(a.2)] Boundary conditions: $m(\emptyset) = 0$ and $\displaystyle\sum_{T \subset G}m(T) = 1$,
                \item [(b.2)] Monotonicity condition: $\forall i \in G$ and $\forall R \subset G\backslash\{i\}$, $m(\{i\}) + \displaystyle\sum_{T \subset R} m(T \cup \{i\}) \geq 0$.
            \end{itemize}
            
            The Choquet integral reformulated in terms of M\"obius representation is defined by
            \begin{equation}\label{choquet_mobius}
                CI(a) = \sum_{T \subseteq G} m(T) min_{j\in T} g_j(a).
            \end{equation}
            
            With the aim of reducing the number of parameters to be elicited, \cite{GRABISCH1997167} introduced the concept of $k$-additive capacity. A capacity is called $k$-additive if $m(T) = 0$ for $T\subseteq G$ such that $|T| > k$. Therefore, for instance, a 2-additive capacity is used to represent situations in which there is interaction between any two criteria but not among three or more criteria.
            
            The value that a 2-additive capacity $\mu$ assigns to a set $S \subseteq G$ can be expressed in terms of M\"obius, and it is defined by
            \begin{equation}\label{mobius_capacity_2-ad}
                \mu(S) = \sum_{j \in S} m(\{j\}) + \sum_{\{j,s\} \subseteq S} m(\{j,s\}), \forall S \subseteq G.
            \end{equation}
            
            With regard to 2-additive capacities, proprieties (a.2) and (b.2) are represented, respectively, by (a.3) and (b.3):
            \begin{itemize}
                \item [(a.3)] Boundary conditions: $m(\emptyset) = 0, \displaystyle\sum_{j \in G} m(\{j\}) + \displaystyle\sum_{\{j,s\} \subseteq G} m(\{j,s\}) = 1$,
                \item [(b.3)] Monotonicity condition:
                $
                    \left\{
                        \begin{array}{l}
                        m(\{j\}) \geq 0, \forall j \in G,\\
                        m(\{j\}) + \displaystyle\sum_{s \in S} m(\{s,j\}) \geq 0, \forall j \in G \hbox{ and } \forall S \subseteq G, S \neq \emptyset.
                    \end{array}\right.
                $
            \end{itemize}
            In the case of 2-additive capacities, the Choquet integral, expressed  in terms of M\"obius representation, is given by
            \begin{equation}
                CI(a) = \sum_{j \in G} m(\{j\}) g_j(a) + \sum_{\{j,s\} \subseteq G} m(\{j,s\}) \min\{g_j(a),g_s(a)\}.\label{choquet_mobius2}
            \end{equation}
        
            Another very useful tool to interpret the Choquet integral is given by the generalized interaction index \citep{Grabisch2000}. For any collation $T$ obtained from the set of criteria criterion $j= 1, \ldots, n$, the interaction index is defined by:
            \begin{equation}
                I(T) = \displaystyle\sum_{K \subset G\backslash T}\left[ \frac{(m-|K| - |T|)!|K|!}{(m-|T|+1)!} \times               \left(\sum_{B \subset T} (-1)^{|T| - |B|}\mu(K \cup B)\right)\right].
            \end{equation}
            \noindent Sets $K$ and $B$ denote auxiliary sets that act as indexes for the sum operators. 
            
            When the generalized interaction index is calculated for a single criterion, say $g_j$, it gives rise to the Shapley importance index, which is denoted by $I_j$ = $I(g_j)$. The Shapley index $I_j$ quantifies the average importance of a single criterion $g_j$ for the decision process. It is always positive and the sum of the Shapley indexes in all criteria is equal to 1, that is, $\sum_{j=1}^{n} I_j = 1$ \citep{Grabisch2000}. For the 2-additive capacities case, the interaction index between two criteria $T=\{g_j, g_s\}$ is represented by $I_{js}$. When $I_{ij}>0$, it represents the synergy between $g_j$ and $g_s$. When representing the redundancy between two criteria $g_j$ and $g_s$, $I_{js}<0$. Finally, $I_{js} = 0$ means that criteria $g_j$ and $g_s$ do not interact.
            
            In 2-additive capacities, monotonicity is ensured when the following restrictions hold:  
            \begin{equation}
                I_{j} - \frac{1}{2}\sum_{j\neq s}|I_{js}| \geq 0, \forall j.
            \end{equation}
            \noindent Moreover, the Choquet integral assumes an insightful form in the case of 2-additive capacities \citep{Grabisch2008}, expressed by
            \begin{equation}
                CI(a) = \sum_{I_{js} > 0} \min \{g_j(a), g_s(a)\}I_{js} + \sum_{I_{js} < 0} \max \{g_j(a), g_s(a)\}|I_{js}| + \sum_{j=1}^{n} g_{j}(a) (I_{j} - \frac{1}{2}\sum_{j\neq s}|I_{js}|).\label{choquet_mobius3}
            \end{equation}
            
            Taking into account the M\"obius representation of a 2-additive capacity, the interaction index is given by $I_{js}=m(\{j,s\})$, and the Shapley value can be written as follows:
            
                $$I_j = m(\{j\}) + \sum_{s \in G\backslash\{j\}}\displaystyle\frac{m(\{j,s\})}{2}.$$
            For a broader discussion of the different representations of the fuzzy measure (capacity, M\"obius and generalized interaction index) see \cite{Grabisch1997}.
            
      \subsection{Commensurability assumption}
      
        \noindent After the introduction of the different Choquet integral formulations in Section \ref{secao_choquet_integral}, we can confirm that a fundamental condition for their application is to express the criteria evaluations on a common scale. This requirement is called commensurability assumption in the fuzzy community \citep{Modave1998}. Indeed, when applying formulation \eqref{eq_choquet}, it is necessary to order the performance evaluations of an alternative, from the smallest to the largest; when applying formulations \eqref{choquet_mobius}, \eqref{choquet_mobius2} and \eqref{choquet_mobius3}, it becomes necessary to identify the minimum value of the performance evaluations for all pairs of interacting criteria to then compute the differences of the criteria evaluations.
        
        Thereby, when all criteria are naturally evaluated on the same scale, the Choquet integral can be applied without any data conversion. On the contrary, the DM has either to assume a qualitative scale to assess the criteria \citep{Modave1998} or to adopt a normalization procedure in order to convert the heterogeneous scales into a common scale \citep{Angilella2015a}. 
    
        To illustrate those different procedures, let us consider a decision-making problem of buying a car, with three purchase options that must be evaluated regarding their maximum speed, in kilometers per hour (km/h), and their consumption, in kilometers per liter (km/l), as presented in Table \ref{tab:exemplo2}.
        
        \begin{table}[htb!]
            \centering
            \small
            \caption{Evaluations of three cars regarding maximum speed and consumption.}
            \begin{tabular}{c|cc}\hline
                   Car    & Maximum Speed (km/h) & Consumption (km/l)\\\hline
                  $a_1$   &  210  &  10  \\
                  $a_2$   &  202  &  12 \\
                  $a_3$   &  200  &  12  \\\hline
            \end{tabular}
            \label{tab:exemplo2}
        \end{table}
        
        Initially, the Choquet integral cannot be applied since these two criteria are evaluated in different scales. Overcoming this limitation, a qualitative scale may be used to assess both criteria, for instance, the scale ``3--good'', ``2--medium'' and ``1--bad''. Evaluations of the cars based on this scale are presented in Table \ref{tab:exemplo1}. Considering the following Shapley and interaction indexes $I_1 = 0.5$, $I_2 = 0.5$ and $I_{12}$ = 0.2, the Choquet integral is given by $CI(a_1) = 2.4 $, $CI(a_2) = 2.4$ and $CI(a_3) = 1.8$. 
        \begin{table}[htb!]
            \centering
            \small
            \caption{Assessing the criteria maximum speed and consumption using a qualitative scale.}
            \begin{tabular}{c|cc}\hline
                   Car    & Maximum Speed (km/h) & Consumption (km/l) \\\hline
                  $a_1$   &  3  &  2  \\
                  $a_2$   &  2  &  3 \\
                  $a_3$   &  1  &  3  \\\hline
            \end{tabular}
            \label{tab:exemplo1}
        \end{table}
        
        
        
        However, as quantitative data is available and gives more detailed information than the qualitative scale, the DM may consider it a loss of information not to use such data in the analysis. Therefore, to keep using the original data, a normalization procedure may be applied in order to convert the heterogeneous scales into a common scale.
        
        Applying, for instance, the min-max normalization ($x' = \frac{x- \min_X}{\max_X - \min_X}$ where $x$ and $x'$ are the original and the normalized values the attribute $X$, respectively, and $\min_X$ and $\max_X$ are the minimum and the maximum values of $X$ \citep{book-han2006}), we acquire the values presented in Table \ref{tab:exemplo3}. The Choquet integral can now be applied, and we obtain $CI(a_1) = 0.4$, $CI(a_2) = 0.52$ and $CI(a_3) = 0.4$. 
        \begin{table}[htb!]
            \centering
            \small
            \caption{Assessing the criteria maximum speed and consumption using a normalized scale.}
            \begin{tabular}{c|cc}\hline
                Cars & Maximum Speed (km/h) & Consumption (km/l)\\\hline
                  $a_1$   &  1        &  0   \\
                  $a_2$   &  0.2      &  1   \\
                  $a_3$   &  0        &  1   \\\hline
            \end{tabular}
            \label{tab:exemplo3}
        \end{table}
        
        \cite{Angilella2015a} proposed a normalization procedure based on simulation and uniform sampling in the range [0,1]: a sample is drawn from a uniform distribution; the generated values are ordered in a decreasing way and these values replace the evaluations of the alternatives respecting the preference order of the DM. For instance, the lowest value generated replace the lowest evaluation, the second lowest value generated replace the second lowest evaluation, and so on and so forth. Then, those generated values (on a common scale) are used on the computation of the Choquet integral instead of the original criteria evaluations. In our example of buying a car, we randomly generated three values for maximum speed and obtained: 0.021, 0.659, 0.524. Consequently, $g_1(a_1) = 0.659$, $g_1(a_2) = 0.524$ and $g_1(a_3) = 0.021$. A drawback of this procedure is to assume that the data is uniformly distributed. It may also be considered not appropriate for some cases since the difference between values impacts the Choquet integral results, and this difference is lost when a uniform sample is generated. 
        
        We can conclude that the choice of the approach to be applied in order to meet the commensurability assumption fully influences the result obtained by the application of the Choquet integral. Therefore, some technical characteristics should be evaluated in order to choose the most appropriate one. For instance, the choice may rely whether the method is compensatory or not, or whether it is an outranking method or a utility-function based method, or even whether a qualitative scale is required, or whether heterogeneous scales can be used without any prior normalization.

    \subsection{The FlowSort sorting method}\label{secao_flowsort}
        
        \noindent FlowSort is a sorting method based on the PROMETHEE methodology for assigning alternatives to $k$ predefined ordered categories $K_1, K_2, \ldots, K_k$, in which $K_1$ is the best category and $K_k$ the worst, i.e, $K_1 \trianglerighteq K_2 \trianglerighteq \ldots \trianglerighteq K_k$. Categories can be defined either by a lower and upper limiting profiles or by central profiles \citep{Nemery200890}. For simplicity, the extension of the FlowSort method proposed in this paper is limited to the case of limiting profiles. 
            
        Let $R = \{r_1, \ldots, r_{k+1}\}$ be the set of reference profiles that characterize the $k$ categories, in which $r_{1}$ and $r_{k+1}$ are the best and the worst reference profiles, respectively. The evaluation of alternatives are also delimited by $r_{1}$ and $r_{k+1}$. Since the categories are completely ordered, each reference profile is preferred to the successive ones, i.e., the conditions 1 and 2 are assumed:
            \begin{flalign*}
                 \hbox{Condition 1:~}& r_1 \succ r_2 \succ \ldots \succ r_k \succ r_{k+1}.\\
                 \mbox{Condition 2:~}&\forall r_{h},~r_{l} \in R, \mbox{ such that }  h < l \Rightarrow g_j(r_{h}) \geq g_j(r_{l}), \forall j= 1, \ldots, n.
            \end{flalign*}
            
        Let us define for any alternative $a_i$ the set $R_i = R \cup \{a_i\}$, $i=1,\ldots,m$. For each criterion, FlowSort builds a preference function $P_j(x,y)=P[d_j(x, y)]$, in which $d_j(x, y) = g_j(x) - g_j(y)$, $x,y \in R_i$. The preference function represents the intensity of preference of $x$ over $y$ on criterion $g_j$, for $j=1, \ldots, n$ and assumes values between 0 and 1. The shape of the preference function for each criterion should be chosen according to the DM preference. Six different types of preference functions are defined by \cite{Brans2005}.
            
            Let us also consider for each criterion $g_j$ a weight $w_j$ that represents its importance, such that $w_j > 0$ and $\sum_{j=1}^n w_j = 1$. For each pair of elements $x,y \in R_i$, the outranking degree, $\pi(x,y)$, is defined as the weighted average of the preference functions as given by:
            \begin{equation}
                \pi(x,y) = \sum_{j=1}^{n}w_j P_j(x,y).\label{outranking_func}
            \end{equation}
            
        The outranking degree $\pi(x,y)$ represents the degree to which element $x$ is preferred to $y$ when considering simultaneously all the criteria. It assumes values between 0 and 1, and the closer $\pi(x,y)$ is to 1, the greater the preference of $x$ over $y$. The outranking degree $\pi(x,y)$, $\forall x, y \in R_i$, follows the conditions 3 to 6:
                \begin{flalign*}
                   \mbox{Condition 3: }& 0 \leq \pi(x,y) \leq 1.\\
                   \mbox{Condition 4: }& \pi(x,y) + \pi(y,x) \leq 1.\\
                   \mbox{Condition 5: }& \pi(x,x)=0.\\
                   \mbox{Condition 6: }& \forall x^{'}, y^{'} \in R_i, \mbox{ if } g_j(x) - g_j(y) \leq g_j(x^{'}) - g_j(y^{'}), \mbox{ then } \pi(x,y) \leq \pi(x^{'}, y^{'}).
                \end{flalign*}
            
            Since the limiting profiles define ordered categories, it can be assumed that if $h <l $, $K_h$ is better than $K_l$. Therefore, profiles must be defined so that $g_j(r_h) > g_j(r_l), ~ \forall j = 1, \ldots, n$ (according
             to Condition 2). Thus, FlowSort imposes that a reference profile of a lower (better) category is preferred to a reference profile of a higher (worse) category. Formally, we have thus Condition 7:
                \begin{flalign*}
                   \mbox{Condition 7: } \forall r_{h},~r_{l} \in R_i, \mbox{ if } h < l \mbox{ then }
                    \pi(r_{h}, r_{l}) > 0 \mbox{ and } \pi(r_{l}, r_{h}) = 0.
                \end{flalign*}
            We can also change a little Condition 7 by accepting that the upper profile $r_h$ of a category $K_h$ is ``strongly preferred'' to the lower profile $r_{h+1}$. This is formalized by Condition 8 and whether reference profiles verify Condition 8, some extra proprieties can be met by FlowSort \citep{thesis_Nemery2008}. 
            \begin{flalign*}
                   \mbox{Condition 8: } \forall r_{h},~r_{l} \in R_i, \mbox{ if } h < l \mbox{ then }
                    \pi(r_{h},r_{l}) = 1.
                \end{flalign*}
            Considering $x, y \in R_i$, positive, negative, and net flows are given by equations \eqref{positive_flow}, \eqref{negative_flow} and \eqref{net_flow}, respectively: 
                \begin{gather}
                    \phi^{+}_{R_i} (x) = \frac{1}{|R_i| - 1}\sum_{\substack{y \in R_i\backslash\{x\}}} \pi(x,y), \label{positive_flow}\\
                    \phi^{-}_{R_i} (x) = \frac{1}{|R_i| - 1}\sum_{\substack{y \in R_i\backslash\{x\} }} \pi(y,x), \label{negative_flow}\\
                    \phi_{R_i} (x) = \phi^{+}_{R_i} (x) - \phi^{-}_{R_i} (x). \label{net_flow}
                \end{gather}
            
            To represent the positive, negative and net flows of an alternative $a_i$, we may simplify the notations for $\phi^{+}(a_i)$, $\phi^{-}(a_i)$ and $\phi(a_i)$, respectively. To assign an alternative $a_i$ to a specific category, its positive and negative flows are compared with the positive and negative flows of the reference profiles, based on the assignment rules given by equations \eqref{comp_positive_flow} and \eqref{comp_negative_flow}, respectively
            \begin{gather}
                \mbox{Assignment rule 1: if } \phi^{+}_{R_i} (r_h) \geq  \phi^{+} (a_i) > \phi^{+}_{R_i} (r_{h+1}), \mbox{ then } K_{\phi^{+}}(a_i) = K_h, \label{comp_positive_flow}\\
                \mbox{Assignment rule 2: if } \phi^{-}_{R_i} (r_h) < \phi^{-} (a_i) \leq \phi^{-}_{R_i} (r_{h+1}), \mbox{ then } K_{\phi^{-}}(a_i) = K_h. \label{comp_negative_flow}
            \end{gather}
            
            \noindent Applying the In order to assign each alternative to exactly one category, the rule based on net flow presented in \eqref{comp_net_flow} can be used:
            \begin{equation}
                \mbox{Assignment rule 3: if } \phi_{R_i} (r_h) \geq \phi(a_i) > \phi_{R_i} (r_{h+1}), \mbox{ then } K_{\phi}(a_i) = K_h. \label{comp_net_flow}
            \end{equation}

\section{The proposed method: FlowSort-Choquet}

    \noindent In this present we propose the FlowSort-Choquet method. The idea of the proposed method is to construct a common scale using preference functions. Moreover, instead of a weighted sum, the Choquet integral is apply to compute the flows. Then, the assignments of alternatives, based on the flows, are conducted using the same sorting rules proposed in FlowSort. This is possible because the flows computed by the Choquet integral keep the same proprieties when they are computed by the weighted sum, as we demonstrate in this section.  
    
    We introduce the FlowSort-Choquet by first presenting its required input data and then presenting the steps for its computation. We also introduces a few fundamental properties that FlowShort-Choquet must satisfy in order to comply with the requirements of a sorting problem with completely ordered categories based on reference profiles.

    \subsection{Problem statement and required input information}
    
        \noindent The MCDA/M problem addressed by the proposed method FlowSort-Choquet refers to sort a set of $m$ alternatives $A=\{a_1, \ldots, a_m\}$ to $k$ predefined ordered categories $K_1, K_2, \ldots, K_k$. The alternatives are evaluated with respect to a set of criteria $G = \{g_1, g_2, \ldots, g_n\}$. Let us consider that the criteria can also be denoted by their index ($G=\{1, 2, \ldots, n\}$). As in FlowSort, categories are defined by a lower and upper limiting profiles following Conditions 1 and 2 (presented in Section 2.3), and we denote the set of limiting profiles by $R = \{r_1, \ldots, r_{k+1}\}$. The limiting profiles have to be defined by the DM, and different profiles may be defined for different criteria. 
        
        We denote the performance evaluations of the alternatives by $g_j(a_i)$, $j=1, \ldots, n$ and $i=1, \ldots, m$, and we assume that heterogeneous scales can be used to assess the alternatives on the different criteria. As established in the FlowSort method, we suppose that the performance of all alternatives in $A$ are between the worst and best limiting profiles. We have thus formally that $\forall a_i \in A, \forall g_j \in G: g_j(r_1)\geq g_{j}(a_i) \geq g_{j}(r_{k+1})$. This is not a hard constraint since the DM is free to define them as she wants.
        
        The criteria of $G$ may have different importance weights, and there may also be interaction between any two criteria. To represent the weights and the interactions, which can be due to synergy or redundancy, the DM has to define the 2-additive capacities ($\mu(\{j\})$ and $\mu(\{s, j\})$, for $s, j=1, \ldots, n, ~ s\neq j$) following conditions (a.3) and (b.3), presented in Section 2.1.
        
        A preference function $P_j$ has to be defined for each criterion $g_j$. Six different types of preference function are proposed in the PROMETHEE method, and they are given in \cite{brans1986}. Depending on the preference function chosen, indifference ($q$) and preference ($p$) thresholds have to be defined. The parameters of the model, including indifference and preference thresholds, reference profiles and capacities as well the performance of alternatives, are assumed to be crisp numbers. 
    
    \subsection{FlowSort-Choquet} \label{secao_computation}
        
        \noindent In order to extend the classical FlowSort method to the Choquet framework, we define for $x,y \in R_i$ the Choquet-outranking degree given by 
        \begin{equation}
            CI_{\pi}(x, y) = \sum_{j \in G} m(\{j\}) P_j(x,y) + \sum_{j,s \in G} m(\{j,s\}) \min \{P_j(x,y) + P_s(x,y)\}.\label{choquet_outeranking_degree}
        \end{equation}
        
        The Choquet-outranking degree can also be expressed in terms of the interaction and Shapley indexes, as given by:
        \begin{eqnarray}
            CI_{\pi}(x, y) &= \displaystyle\sum_{I_{js} > 0} \min \{P_j(x,y), P_s(x,y)\}I_{js} + \displaystyle\sum_{I_{js} < 0} \max \{P_j(x,y), P_s(x,y)\} |I_{js}| \nonumber\\ &+ \displaystyle\sum_{j=1}^{n} P_j(x,y) (I_{j} - \displaystyle\frac{1}{2}\displaystyle\sum_{j\neq s}|I_{js}|)
            \label{choquet_outeranking_degree2}
        \end{eqnarray}
        
        Our proposal here is to apply the Choquet-outranking degree $CI_{\pi}(x, y)$ to calculate the flows instead of applying the outranking degree $\pi(x,y)$ defined in \eqref{outranking_func}. That is possible because the Choquet-outranking degree meets conditions 3.B to 7.B:
         \begin{flalign*}
               \mbox{Condition 3.B: }& 0 \leq CI_{\pi}(x, y) \leq 1. \\
               \mbox{Condition 4.B: }& CI_{\pi}(x, y) + CI_{\pi}(y, x)  \leq 1.\\
               \mbox{Condition 5.B: }& CI_{\pi}(x, x) =0.\\
               \mbox{Condition 6.B: }& \forall x^{'}, y^{'} \in R_i, \mbox{ se } g_j(x) - g_j(y) \leq g_j(x^{'}) - g_j(y^{'}), \mbox{ then } CI_{\pi}(x, y)  \leq CI_{\pi}(x', y') .\\
               \mbox{Condition 7.B: }& \forall r_{h},~r_{l} \in R_i, \mbox{ if } h < l \mbox{ then }
                    CI_{\pi}(r_{h}, r_{l}) > 0 \mbox{ and } CI_{\pi}(r_{l}, r_{h}) = 0.\\
        \end{flalign*}
        
        \noindent Proof: Indeed, as $0 \leq P_j(x,y) \leq 1$ and by the proprieties (a.3) and (b.3), $CI_{\pi}(x, y)$ meets Condition 3.B. By propriety (a.3) and the fact that if $P_j(x,y) > 0$ then $P_j(y,x) =0$, Condition 4.B is verified. Condition 5.B is verified directly from the fact that $P(x,x) = 0$. Condition 6.B is met since $g_j(x) - g_j(y) \leq g_j(x^{'}) - g_j(y^{'})$ implies $P(x,y) \leq P(x',y')$, and that, along with proprieties (a.3) and (b.3), point to $CI_{\pi}(x, y)) \leq CI_{\pi}(x', y')$. As a reference profile of a lower (better) category is preferred to a reference profile of a higher (worse) category, we have Condition 7.B. 
        \begin{flushright}
        $\square$
        \end{flushright}
        
        As in FlowSort, we can change Condition 7.B by accepting that the upper profile $r_h$ of a category $K_h$ is ``strongly preferred'' to the lower profile $r_{h+1}$, resulting in Condition 8.B. If reference profiles verify Condition 8.B, then some extra proprieties can be met by Choquet-FlowSort, as presented in Section \ref{section_proprieties}.
        \begin{flalign*}            
            \mbox{Condition 8.B: }& \forall r_{h},~r_{l} \in R_i, \mbox{ if } h < l \mbox{ then }
            CI_{\pi}(r_{h},r_{l}) = 1.
        \end{flalign*}

        Therefore, considering the capacities $\mu$ instead of the weights $w$, and using the Choquet-outranking degree as defined in \eqref{choquet_outeranking_degree} or \eqref{choquet_outeranking_degree2}, the positive and negative flows in the FlowSort-Choquet, called Choquet-flows, are given, respectively, by equations
        \begin{eqnarray}
            \phi^{+}_{R_i,CI_{\pi}}(x) = \frac{1}{|R_i| - 1}\sum_{\substack{y \in R_i\backslash\{x\}}} CI_{\pi}(x, y),
            \label{pos_flow_flow-choquet2} \hbox{ and }\\
            \phi^{-}_{R_i,CI_{\pi}}(x) = \frac{1}{|R_i| - 1}\sum_{\substack{y \in R_i\backslash\{x\}}} CI_{\pi}(y, x)).\label{neg_flow_flow-choquet2}
        \end{eqnarray}
        
        \noindent The net Choquet-flow is then given by
        \begin{equation}
            \phi_{R_i,CI_{\pi}}(x) = \phi^{+}_{R_i,CI_{\pi}}(x) - \phi^{-}_{R_i,CI_{\pi}}(x).\label{net_flow_flow-choquet2}
        \end{equation}
        
        \begin{propos} \label{prop_1} Under conditions 1 to 6, we have that the order of the Choquet-flows of the reference profiles is invariant with respect to alternative $a_i \in A$, i.e, $\forall a_i \in A$ and $\forall h \in \{1, \ldots, k\}$
            \begin{eqnarray}
                \phi_{R_i,CI_{\pi}}^{+} (r_h) > \phi_{R_i,CI_{\pi}}^{+} (r_{h+1}) \label{prova_1}\\
                \phi_{R_i,CI_{\pi}}^{-} (r_h) < \phi_{R_i,CI_{\pi}}^{-} (r_{h+1}) \\
                \phi_{R_i,CI_{\pi}} (r_h) > \phi_{R_i,CI_{\pi}} (r_{h+1})
            \end{eqnarray}
        \end{propos}
        Proof: By condition 2, $g_j(r_h)\geq g_j(r_{h+1}), \forall j \in G, \forall h=1, \ldots, k.$ Consequently, by condition 6.B, 
        \begin{equation} CI_{\pi}(r_h, r_l) \geq CI_{\pi}(r_{h+1}, r_l).\label{eq_1}\end{equation}
        Conditions 7.B, $CI_{\pi}(r_h, r_{h+1}) > 0$, and 5.B, $CI_{\pi}(r_{h+1}, r_{h+1}) = 0$, imply
        \begin{equation}
            CI_{\pi}(r_h, r_{h+1}) > CI_{\pi}(r_{h+1}, r_{h+1}).\label{eq_2}
        \end{equation}
        By equations \eqref{eq_1} and \eqref{eq_2} we have
        \begin{equation}\label{eq_3}
            \phi_{R,CI_{\pi}}^{+} (r_h) = \sum_{r_l \in R} CI_{\pi}(r_h, r_l) \geq \sum_{r_l \in R} CI_{\pi}(r_{h+1}, r_l) =  \phi_{R,CI_{\pi}}^{+} (r_{h+1}).
        \end{equation}
        Combining Condition 1 with Condition 6.B implies that: $\forall a_i \in A, CI_{\pi}(r_h, a_i) \geq CI_{\pi}(r_{h+1}, a_i)$. By conditions 2 and 6.B we have that $CI_{\pi}(r_h, a_i) \geq CI_{\pi}(r_{h+1}, a_i)$. Since 
        \begin{flalign*}
            &\phi_{R_i,CI_{\pi}}^{+}(r_h)= \phi_{R,CI_{\pi}}^{+}(r_h) + CI_{\pi}(r_h, a_i),\\
            &\phi_{R_i,CI_{\pi}}^{+}(r_{h+1})= \phi_{R,CI_{\pi}}^{+}(r_{h+1}) + CI_{\pi}(r_{h+1}, a_i), \hbox{~and }\\
            &\phi_{R ,CI_{\pi}}^{+}(r_h) > \phi_{R,CI_{\pi}}^{+}(r_{h+1}), 
        \end{flalign*}
        we can conclude that $\phi_{R_i ,CI_{\pi}}^{+}(r_h) > \phi_{R_i,CI_{\pi}}^{+}(r_{h+1})$ and expression \eqref{prova_1} is proven. The proofs for the negative and net Choquet-flows are similarly obtained.
        \begin{flushright}$\square$\end{flushright}
        
        Therefore, by Proposition \ref{prop_1}, although the Choquet-flows of the reference profiles directly depend on the alternative $a_i$, their order always respects the order of the categories. This allows us to delimit a category $K_h$ by the Choquet-flow values of $r_h$ and $r_{h+1}$ in the case that the categories are defined by an upper and lower limit. Therefore, as in FlowSort \citep{Nemery200890}, this proposition is the basis of the Choquet-FlowSort assignment rules. We define here three assignment rules:
        \begin{eqnarray}
            &\hbox{Assignment rule 1:}& \hbox{if } \phi_{R_i,CI_{\pi}}^{+} (r_h) \geq \phi_{CI_{\pi}}^{+} (a_i) > \phi_{R_i,CI_{\pi}}^{+} (r_{h+1}), \hbox{~then~} K_{\phi^{+}} (a_i) = K_{h}, \label{rule1}\\
            &\hbox{Assignment rule 2:}& \hbox{if } \phi_{R_i,CI_{\pi}}^{-} (r_h) < \phi_{CI_{\pi}}^{-} (a_i) \leq \phi_{R_i,CI_{\pi}}^{-} (r_{h+1}), \hbox{~then~} K_{\phi^{-}} (a_i) = K_{h}, \label{rule2}\\
            &\hbox{Assignment rule 3:}& \hbox{if } \phi_{R_i,CI_{\pi}} (r_h) \geq \phi_{CI_{\pi}} (a_i) > \phi_{R_i,CI_{\pi}} (r_{h+1}), \hbox{~then~} K_{\phi} (a_i) = K_{h}. \label{rule3}
        \end{eqnarray}
        
\begin{figure}[h!]
             \centering
             \includegraphics[width=12cm]{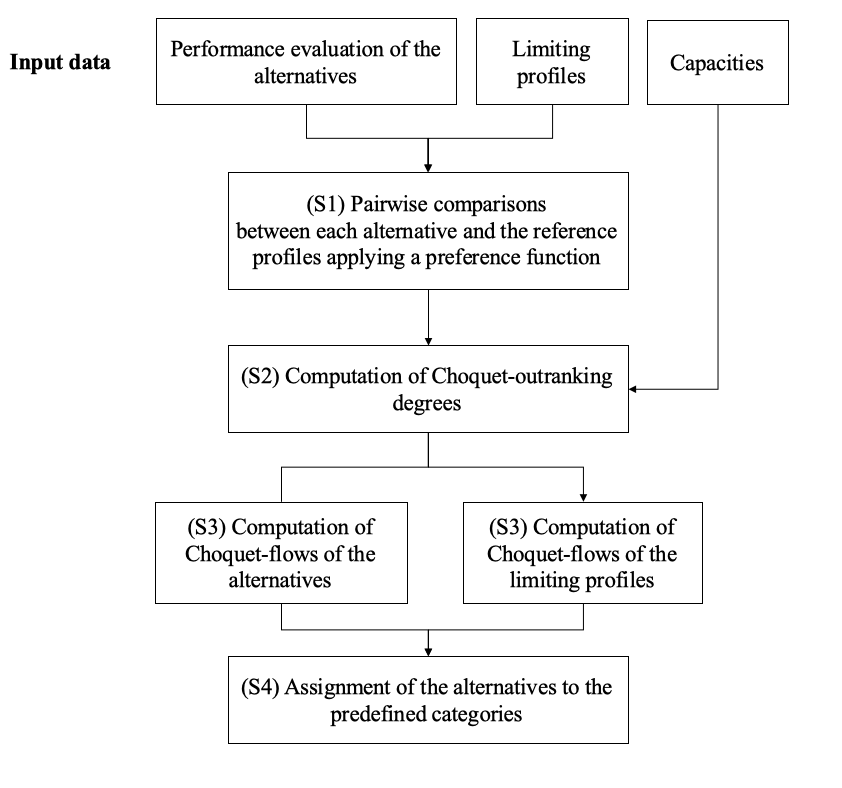}
             \caption{The FlowSort-Choquet algorithm scheme.}
             \label{fig:Algoritmo}
         \end{figure}

        After all, we can introduce the steps of FlowSort-Choquet, as illustrated in Figure \ref{fig:Algoritmo}. Since all required input data has been defined, the first step (S1) of the FlowSort-Choquet algorithm is to conduct the pairwise comparisons between each alternative and the reference profiles applying the preference function $P_j(x,y)=P[d_j(x, y)]$, in which $d_j(x, y) = g_j(x) - g_j(y)$, $x,y \in R_i$ and $R_i = R \cup \{a_i\}$, $i=1,\ldots, m$. In this step, we obtain a common scale, since $P_j(x,y)$ assumes values between 0 and 1, from a heterogeneous scale (for the different criteria $g_j \in G$). The second step is the computation of the Choquet-outranking degrees using equation \eqref{choquet_outeranking_degree} or \eqref{choquet_outeranking_degree2} for all elements $x, y \in R_i, ~i=1, \ldots, m$. 
         
        In the third step (S3), the positive and negative Choquet-flows of the alternatives and limiting profiles are computed  by using equations \eqref{pos_flow_flow-choquet2} and \eqref{neg_flow_flow-choquet2}, respectively, and the Choquet-net flow by equation \eqref{net_flow_flow-choquet2}. In the fourth step (S4), the assignments may be made using the Choquet-flows and the assignment rules defined in \eqref{rule1}, \eqref{rule2} and \eqref{rule1}. 
         
    \section{Analysis of some FlowSort-Choquet proprieties}\label{section_proprieties}
    
        In order to comply with the requirements of a sorting problem with completely ordered categories based on reference profiles, the FlowShort-Choquet method must satisfy certain fundamental properties, which are introduced in this Section. Referring to the properties of Yu's ordinal method \citep{Yu1992}, five fundamental principles  have been defined which characterize a sorting method: principles of universality, independency, neutrality, stability and homogeneity. The fulfillment of these principles is also discussed by \cite{thesis_Belacel2000} for the PROAFTN sorting method and by \cite{thesis_Nemery2008} for FlowSort.

        The principle of universality says that the sorting procedure may assign an alternative to one, many, or no categories. Particularly, the FlowSort-Choquet method assigns each alternative to exactly one category, following the principle of uniqueness (Proposition \ref{prop_uniqueness}), which is a particular case of the universality principle. The principles of independency and neutrality are presented in Propositions \ref{propos_independency} and \ref{prop_neutrality}, respectively. 

        \begin{propos}\label{prop_uniqueness}
        \textbf{\emph{Propriety of uniqueness}}\\
            The sorting method assigns each alternative $a_i$ to exactly one category.
        \end{propos}
        Proof: According to each assignment rule \eqref{rule1}, \eqref{rule2}, \eqref{rule3}, each alternative $a_i \in A$ is assigned to exactly one category. The property of uniqueness is thus fulfilled.
        \begin{flushright}
        $\square$
        \end{flushright}

        \begin{propos}\label{propos_independency}
        \textbf{\emph{Propriety of independency}}\\
            The assignment of an alternative $a_i$ does not depend on the assignment of another alternative $a_t$.
        \end{propos}
        Proof: FlowSort-Choquet assigns an alternative $a_i \in A$ to a category by comparing it only with the reference profiles and regardless of the other alternatives of $A$. Hence, the principle independence is verified.
        \begin{flushright}
        $\square$
        \end{flushright}
        
        \begin{propos}\label{prop_neutrality}
        \textbf{\emph{Propriety of neutrality}}\\
            The assignment of each alternative does not depend on its given label, i.e, if we give two different labels
            $a_1$ and $a_2$ to the same alternative, FlowSort-Choquet assigns both of them to the same category:
            $$K_{\phi}(a_1) = K_{\phi}(a_2).$$
        \end{propos}
        Proof: It is verified by Proposition \ref{propos_independency} and by the definition of the assignment rules themselves. 
        \begin{flushright}
        $\square$
        \end{flushright}
        
        In the sorting context, the principle of stability is usually defined as the fact that the assignment of an alternative is not altered when one or several categories (different from the category to which the alternative is assigned) are suppressed or added. As discussed by \cite{thesis_Nemery2008}, since the assignments of the alternatives depend on the definition of the categories (particularly, on the definition of the reference profiles), this property is not respected in FlowSort. Nevertheless, FlowSort assignment rules respect what is called weak stability. 
        
        Under the same justification, FlowSort-Choquet respects only the weak stability propriety, formally defined in Proposition \ref{prop_stability-1}, in case of splitting of the categories, and in Proposition \ref{prop_stability-1B} in case of fusion of the categories. We present those propositions considering the positive Choquet-Flows, but for the negative and net Choquet-flows them can be similarly verified. As the Choquet-outranking degree $CI_{\pi}(x,y)$, $\forall x, y \in R_i, \forall i=1, \ldots, m$, complies with all the same conditions required by the outranking degree $\pi(x,y)$, the proofs are analogous to those presented by \cite{thesis_Nemery2008} for the FlowSort method. 
        
        \begin{propos}\textbf{\emph{Propriety of weak stability-1: splitting of categories}}\label{prop_stability-1}\\
        Let $a_i$ be any of the alternatives of $A$, and let us suppose that $K_{\phi^{+}} (a_i) = K_{h}$. Under Condition 8.B we have: 
        \begin{eqnarray*}
            \hbox{~If ~} R^{'} = R \backslash \{r_s\} \hbox{~with~} s < h \hbox{~and~} K_{\phi^{+}, R^{'}} (a_i) = K^{'}_{h^{'}} \implies h^{'} \in [h-1, h].\\
            \hbox{~If ~} R^{'} = R \backslash \{r_s\} \hbox{~with~} h < s \hbox{~and~} K_{\phi^{+}, R^{'}} (a_i) = K^{'}_{h^{'}} \implies h^{'} \in [h+1, h].
        \end{eqnarray*}
        \end{propos}
        Proof: Analogous to the demonstration presented by \cite{thesis_Nemery2008}.
        \begin{flushright}
        $\square$
        \end{flushright}
        
        \begin{propos}\textbf{\emph{Propriety of weak stability-1: fusion of categories}}\label{prop_stability-1B}\\
        Let $a_i$ be any of the alternatives of $A$, and let us suppose that $K_{\phi^{+}} (a_i) = K_{h}$. Under Condition 8.B we have: 
        \begin{eqnarray*}
            \hbox{~If ~} R^{'} = R \cup \{r_s\} \hbox{~with~} CI_\pi(r^{'}_s, r_h) = 1 \hbox{~and~} K_{\phi^{+}, R^{'}} (a_i) = K^{'}_{h^{'}} \implies h^{'} \in [h, h+1].\\
            \hbox{~If ~} R^{'} = R \cup \{r_s\} \hbox{~with~} CI_\pi(r_h, r^{'}_s) = 1 \hbox{~and~} K_{\phi^{+}, R^{'}} (a_i) = K^{'}_{h^{'}} \implies h^{'} \in [h-1, h].
        \end{eqnarray*}
        \end{propos}
        Proof: Analogous to the demonstration presented by \cite{thesis_Nemery2008}.
        \begin{flushright}
        $\square$
        \end{flushright}
        
        The homogeneity principle says that if the outranking (preference) relations between an alternative $a_i$ and the reference profiles are the same as the outranking relations between $a_t$ and the reference profiles, then $a_i$ and $a_t$ are assigned to the same categories. In the FlowSort-Choquet method, as in FlowSort, the weak-homogeneity is fulfilled in terms of preference degrees (and not preference relations), as showed in Proposition \ref{prop_homo}.
    
        \begin{propos}\textbf{\emph{Property of weak homogeneity}}\label{prop_homo}\\
        If the preference (outranking) degrees between an alternative $a_i$ and the reference profiles are the same as the preference degrees between $a_t$ and the reference profiles, FlowSort-Choquet assigns $a_i$ and $a_t$ to the same category, i.e, $\forall a_i, a_t \in A$ and $\forall r_h \in R$,
        \begin{eqnarray*}
        &\hbox{if}&  CI_{\pi}(r_h, a_i) = CI_{\pi}(r_h, a_t) \hbox{ and } CI_{\pi}(a_i, r_h) = CI_{\pi}(a_t, r_h) \\
        &\implies& K_{\phi^{+}} (a_i) = K_{\phi^{+}}(a_t), ~ K_{\phi^{-}} (a_i) = K_{\phi^{-}} (a_t) \hbox{ and } K_{\phi} (a_i) = K_{\phi} (a_t).
        \end{eqnarray*}
        \end{propos}
        Proof: We can see that the computation of the flows (equations \eqref{pos_flow_flow-choquet2}, \eqref{neg_flow_flow-choquet2}, \eqref{net_flow_flow-choquet2}) depends directly on the Choquet-outranking degree defined in equation \eqref{choquet_outeranking_degree}. Since the Choquet-outranking degrees are equal ($CI_{\pi}(r_h, a_i) = CI_{\pi}(r_h, P(r_h,a_t)) \hbox{ and } CI_{\pi}(a_i, P(a_i,r_h)) = CI_\pi(a_t, P(a_t,r_h))$), we have that the Choquet-flows of any two alternatives $a_i$ and $a_j$, with respect to the reference profiles, will thus be equal, i.e, $\phi_{R_i,CI_{\pi}}^{+} (a_i) = \phi_{t, CI_\pi}^{+} (a_t)$, $\phi_{R_i,CI_\pi}^{-} (a_i) = \phi_{R_t,CI_{\pi}}^{-} (a_t)$ and $\phi_{i, CI_\pi} (a_i) = \phi_{R_t,CI_{\pi}} (a_t)$. Similarly, the flows taken by the profiles with respect to the two alternatives are also equal: $\forall r_h \in R$, $\phi_{R_i,CI_{\pi}}^{+} (r_h) = \phi_{R_t,CI_{\pi}}^{+} (r_h)$, $\phi_{R_i, CI_\pi}^{-} (r_h) = \phi_{R_t,CI_\pi}^{-} (r_h)$ and $\phi_{R_i,CI_{\pi}} (r_h) = \phi_{R_t,CI_{\pi}} (r_h)$. Consequently, by the assignments rules \eqref{rule1}, \eqref{rule2} and \eqref{rule3}, $K_{\phi^{+}} (a_i) = K_{\phi^{+}}(a_t), ~ K_{\phi^{-}} (a_i) = K_{\phi^{-}} (a_t) \hbox{ and } K_{\phi} (a_i) = K_{\phi} (a_t)$. 
        \begin{flushright}
        $\square$
        \end{flushright}
        
        Another four principles (principle of pairwise assignment consistency, stability-2, monotonicity and category conformity) have been defined to characterize sorting procedures with completely ordered categories based on limiting profiles as presented by \cite{thesis_Nemery2008}. 
        
        A sorting procedure fulfill the pairwise assignment consistency propriety when its preferential relational system respects pairwise comparisons between two alternatives regardless of the sorting problem, i.e, if preference relations are transitive. In a sorting procedure, we say that the preference relation is transitive if the assignment of $a_i$ to a better category than $a_t$ implies that $a_i$ is preferred to $a_t$. However, due to the Paradox of Condorcet \citep{thesis_Nemery2008, thesis_Belacel2000}, in the FlowSort-Choquet method, as in outranking methods in general, preference relations are not necessarily transitive and an alternative $a_i$ may be assigned to a worse category than $a_t$ although we have that $a_t$ is preferred to $a_i$. 
        
        This situation may be illustrated with a decision-problem with two categories ($K_1$, $K_2$) and with the following alternatives and limiting profile evaluated on 3 criteria (to be maximized and with same weights): $g_1(a_i)= 3, ~g_2(a_i)=3, ~g_3(a_i)=1$, and $g_1(a_t)=4, ~g_2(a_t)=1, ~g_3(a_t)=2$, and $g_1(r_2)=1, g_2(r_2)=2, g_3(r_2)=3$. Since we have $a_i \succ r_2$ and $r_2 \succ a_t$, the assignments are as follows: $K(a_i) = K_1$ and $K(a_t) = K_2$. Nevertheless, we have that $a_t \succ a_i$. Consequently, analog to FlowSort, due to the non-transitivity of the preference and outranking relations, Flowsort-Choquet does not fulfill the pairwise assignment consistency propriety either. 
        
        Regarding the principle of stability-2, its definition is similar to the propriety of weak stability-1 (Propositions \ref{prop_stability-1} and \ref{prop_stability-1B}), except by the fact that the fusion or the splitting is limited to consecutive categories. The FlowSort-Choquet meets the weak stability-2 proprieties under some conditions, as presented in Propositions \ref{prop_stability-2} and \ref{prop_stability-2B}.
        
        \begin{propos}\textbf{\emph{Propriety of weak stability-2: splitting of categories}}\label{prop_stability-2}\\
        Under Condition 8.B, the assignment of an alternative $a_i$ to $C_h$, according to the positive Choquet-flows, will not be affected by the splitting of two consecutive categories if the following conditions are satisfied:
        \begin{eqnarray}
            \phi_{CI_{\pi}}^{+} (a_i) - \phi_{R_i, CI_{\pi}}^{+} (r_{h+1}) &>& 
            \displaystyle\frac{CI_{\pi}(r_{h+1}, r_l)}{k+1} - \displaystyle\frac{CI_{\pi}(a_i, r_l)}{k+1},\\
            \phi_{R_i, CI_{\pi}}^{+} (r_h) - \phi_{CI_{\pi}}^{+} (a_i) &>& 
            \displaystyle\frac{CI_{\pi}(a_i, r_l)}{k+1} - \displaystyle\frac{CI_{\pi}(r_{h}, r_l)}{k+1}.
        \end{eqnarray}
        \end{propos}
        Proof: Analogous to the demonstration presented by \cite{thesis_Nemery2008}.
        \begin{flushright}
        $\square$
        \end{flushright}
        
        \begin{propos}\textbf{\emph{Propriety of weak stability-2: fusion of categories}}\label{prop_stability-2B}\\
        Under Condition 8.B, the assignment of an alternative $a_i$ to $C_h$, according to the positive Choquet-flows, will not be affected by the fusion of two consecutive categories if the following conditions are satisfied:
        \begin{eqnarray}
            \phi_{R_i,CI_{\pi}}^{+} (r_h) - \phi_{CI_{\pi}}^{+} (a_i) &>& 
            \displaystyle\frac{CI_{\pi}(r_h, r_l)}{k+1} - \displaystyle\frac{CI_{\pi}(a_i, r_l)}{k+1},\\
            \phi_{CI_{\pi}}^{+} (a_i) - \phi_{R_i, CI_{\pi}}^{+} (r_{h+1}) &>& 
            \displaystyle\frac{CI_{\pi}(a_i, r_l)}{k+1} - \displaystyle\frac{CI_{\pi}(r_{h+1}, r_l)}{k+1}.
        \end{eqnarray}
        \end{propos}
        Proof: Analogous of the demonstration presented by \cite{thesis_Nemery2008}.
        \begin{flushright}
        $\square$
        \end{flushright}
        
        The two other principles, monotonicity and category conformity, are completely fulfilled by FlowSort-Choquet as presented in Propositions \ref{prop_monotonicity} and \ref{prop_categoryconformity}, respectively. 
     
        \begin{propos}\textbf{\emph{Propriety of monotonicity}}\label{prop_monotonicity}\\
        If  an  alternative $a_i$ dominates another alternative $a_t$,  then $a_i$ can  not be assigned to a higher (worse) category than alternative $a_t$. Formally, $\forall a_i, a_t \in A$, 
        \begin{eqnarray*}
               &\hbox{if}& g_j(a_i) \geq g_j(a_t)  ~\forall g_j \in G \hbox{ and } \exists g_x(a_t) \leq g_x(a_i) \\&\implies& K_{\phi^{+}} (a_i) \trianglerighteq K_{\phi^{+}} (a_t) ~\hbox{and} ~ K_{\phi^{-}} (a_i) \trianglerighteq K_{\phi^{-}} (a_t) ~\hbox{and}~K_{\phi} (a_i) \trianglerighteq K_{\phi} (a_t)
            \end{eqnarray*}
        \end{propos}
        Proof: Let us suppose that $a_t$ is assigned to category $K_h$. In such a case, by the assignment rule \eqref{rule1}, we must have that 
        \begin{equation}\label{prova1}
            \phi_{R_t, CI_{\pi}}^{+} (a_t) > \phi_{R_t, CI_{\pi}}^{+}(r_{h+1}). 
        \end{equation}
        Conversely, since $g_j(a_i) \geq g_j(a_t), \forall j=1\ldots, n$, we have that $-g_j(a_i) \leq -g_j(a_t) \implies g_j(r_{h+1}) - g_j(a_i) \leq g_j(r_{h+1}) - g_j(a_t)$ and by Condition 4.B we have that 
        \begin{eqnarray}\label{prova2}
            &~& CI_{\pi}(r_{h+1},a_t) \leq CI_{\pi}(r_{h+1},a_t)) \nonumber\\
            &\implies &(|R_i| - 1)\phi_{R}^{+} (r_{h+1}) + CI_{\pi}(r_{h+1}, a_t)) \leq (|R_i| - 1)\phi_{R}^{+} (r_{h+1}) + CI_{\pi}(r_{h+1}, a_t) \nonumber\\
            &\implies & \phi_{R_i}^{+} (r_{h+1}) \leq \phi_{R_t}^{+} (r_{h+1}) \label{prova2}
        \end{eqnarray}
        Similarly, as $g_j(a_i) \geq g_j(a_t), \forall j=1\ldots, n$, then $g_j(a_i) - g_j(r) \geq g_j(a_t) - g_j(r), \forall r \in R$, and by Condition 6.B, $CI_{\pi}(a_i, r)) \geq CI_{\pi}(a_t, r)), \forall r \in R.$ Since by Condition 5.B $CI_{\pi}(a_i, a_i)) = 0$ and $CI_{\pi}(a_t, a_t)) = 0$, we obtain that  
        \begin{eqnarray}
            \sum_{r \in R_i} CI_{\pi}(a_i, P(a_i, r)) \geq \sum_{r \in R_t} CI_{\pi}(a_t, r)
            \implies ~ \phi_{R_i, CI_{\pi}}^{+} (a_i) \geq \phi_{R_t,CI_{\pi}}^{+} (a_t). \label{prova3}
        \end{eqnarray}
           From expressions \eqref{prova3}, \eqref{prova1} and \eqref{prova2}, respectively, we obtain
        $$\phi_{R_i,CI_{\pi}}^{+} (a_i) \geq \phi_{R_t,CI_{\pi}}^{+} (a_t) \geq \phi_{R_t,CI_{\pi}}^{+} (r_{h+1}) \geq \phi_{R_i,CI_{\pi}}^{+} (r_{h+1}).$$
         Consequently, 
        $$K_{\phi^{+}}(a_i) \trianglerighteq K_h = K_{\phi^{+}}(a_t).$$
        This proves the proposition considering the positive Choquet-flow assignment. The proof for the negative and net Choquet-flows assignments are analogous.
        \begin{flushright}
        $\square$
        \end{flushright}
        
        \begin{propos}\textbf{\emph{Property of category conformity}}\label{prop_categoryconformity}\\
        If the performances of an alternative $a_i$ (on all the criteria) are in between two successive limiting profiles, then the alternative $a_i$ may be assigned to the corresponding category, i.e, $\forall a_i \in A$ and $\forall g_j \in G$, 
        $$\hbox{if } g_j(r_{h+1}) \leq g_j(a_i) \leq g_j(r_h), \hbox{then } K_{\phi^{+}} (a_i) = K_{\phi^{-}} (a_i) = K_{\phi} (a_i) = K_h.$$
        \end{propos}
        Proof: The proof is immediate from proposition \ref{prop_1}.
        \begin{flushright}
        $\square$
        \end{flushright}
        
        With the aim of explaining the assignments of the alternatives in different situations, we present two more FlowSort-Choquet proprieties related to the assignment rules of FlowSort-Choquet, regarding the relationship between $K_{\phi^+}$ and $K_{\phi^-}$ and the coherence of the net-flow assignment rule, which are proved in a similar way to FlowSort \citep{thesis_Nemery2008}.
    
        \begin{propos} \textbf{\emph{Relationship between $K_{\phi^+}$ and $K_{\phi^-}$}}\\
        Under Condition 8.B we have that category $K_{\phi^{-}}(a_i)$ is always as least as good as category $K_{\phi^{+}}(a_i)$ .
        \end{propos}
        Proof: As the Choquet-outranking degree $CI_{\pi}(x,y)$, $\forall x, y \in R_i, \forall i=1, \ldots, m$, complies with all the same conditions required by the preference degree $\pi(x,y)$, the proof is analogous to the one presented in \cite{thesis_Nemery2008} for the FlowSort method. 
        \begin{flushright}
        $\square$
        \end{flushright}
        
        \begin{propos}
            \textbf{\emph{Coherence of the net-flow assignment rule}}\\
            The assignment procedure based on the net Choquet-flows complies with the one based on the positive and negative Choquet-flows since the assignment result is always in between the positive and the negative assignment result. More formally, we have, $\forall a_i \in A$
            \begin{eqnarray*}
                K_{\phi^{-}} (a_i) \trianglerighteq K_{\phi} (a_i) \trianglerighteq K_{\phi^{+}} (a_i).
            \end{eqnarray*}
        \end{propos}
        Proof: As the Choquet-outranking degree $CI_{\pi}(x,y)$, $\forall x, y \in R_i, \forall i=1, \ldots, m$ complies with all the same conditions required by the preference degree $\pi(x,y)$, the proof is analogous to the one presented in \cite{thesis_Nemery2008} for the FlowSort method. 
        \begin{flushright}
        $\square$
        \end{flushright}
    
    \section{Numerical example} \label{sec_example}
    
    \noindent In order to exemplify the proposed method, we present a decision-making example with interacting criteria, and in which criteria evaluations are expressed in heterogeneous scales. This problem is based on the example presented by \cite{Angilella2015a}.

    Ten cars (the alternatives) are evaluated regarding four criteria: price (in Euro), acceleration (0 to 100 km/h in seconds, that is, how many seconds is necessary to achieve a speed of 100km/h), maximum speed (in km/h) and consumption (in l/100 km). The performance evaluation of the cars are presented in Table \ref{tab_exemplo2_evaluations}. Criteria price, acceleration and consumption have a decreasing direction of preference (the lower, the better), while criterion maximum speed has an increasing direction of preference (the higher, the better).
        
        \begin{table}[htb!]
            \centering
            \small
            \caption{Performance evaluation of the ten car models regarding the four considered criteria: price, acceleration, maximum speed and consumption. }\label{tab_exemplo2_evaluations}
            \begin{tabular}{lcccc}\hline
                	Cars	&	Price	&	Acceleration	&	Max Speed	&	Consumption	\\
                    	    &	(Euro)	&	(sec.0 to 100)	&	(km/h)	&	(l/100 km)	\\\hline
                    $a_1$   &	16,000	&	12.0	&	185	&	3.1	\\
                    $a_2$   &	15,750	&	13.5	&	163	&	3.8	\\
                    $a_3$   &	15,050	&	11.0	&	173	&	4.0	\\
                    $a_4$   &	15,260	&	12.0	&	172	&	3.3	\\
                    $a_5$   &	16,300	&	10.6	&	183	&	3.7	\\
                    $a_6$   &	16,050	&	10.8	&	180	&	3.4	\\
                    $a_7$   &	17,000	&	11.0	&	170	&	3.8	\\
                    $a_8$   &	17,500	&	12.9	&	174	&	3.5	\\
                    $a_9$   &	17,800	&	11.8	&	170	&	3.8	\\
                    $a_{10}$&	17,060	&	13.9	&	175	&	3.9	\\\hline
            \end{tabular}
        \end{table}
     Performing the weights of the criteria, the DM defines the Shapley values $I_1=0.25$, $I_2=0.21$, $I_3=0.16$ and $I_2=0.38$, indicating a greater importance of the criterion consumption, followed by criteria price, acceleration and lastly, with least importance, criterion maximum speed. Moreover, let us suppose that the DM considers that criteria acceleration and maximum speed are redundant, since, in general, speedy cars also have a good acceleration, and that there is a synergy between criteria maximum speed and consumption, since a speedy car also with a low consumption is very well appreciated. To represent those interactions, the DM defines the interaction indexes $I_{23} = -0.08$ and $I_{34} = 0.10$.
        
    The goal here is to obtain a purchase viability assessment of the cars, classifying them into three categories: ($K_1$) very feasible, ($K_2$) feasible and ($K_3$) not feasible. To define those categories, the DM determines four limiting profiles for each criteria as presented in Table \ref{tab:ex2_profiles}.
        
        \begin{table}[htb!]
            \centering\small
            \caption{Limiting profiles defining the categories very feasible ($K_1$), feasible ($K_2$) and not feasible ($K_3$).}
            \label{tab:ex2_profiles}
            \begin{tabular}{l|cccc}\hline
                Limiting Profiles	&	Price	&	Acceleration	&	Max Speed	&	Consumption	\\
                    	&	(Euro)	&	(in sec 0/100 km/h)	&	(km/h)	&	(l/100 km)	\\\hline
                $r_1$     &  15,000         &   10.5                &       190    &       3.0     \\
                $r_2$     &  16,000         &   11.5                &       180    &       3.4        \\
                $r_3$     &  17,000         &   12.5               &        170   &        3.8       \\
                $r_4$     &  18,000         &   15.0                &       160    &       4.2       \\\hline
            \end{tabular}
        \end{table}
        
        The type of the preference function has also to be defined. We consider here the preference function of type 1 (the usual function), given by equation \eqref{eq_pref_function}, for all criteria; i.e, two elements $x$ and $y$ are indifferent if $g_j(x) = g_j(y)$. In that case, the DM does not have to define values for the thresholds. 
                \begin{equation}\label{eq_pref_function}
                    P_j(x,y) = P (d_j (x, y))
                        \left\{
                        	\begin{array}{ll}
                        		0, & \mbox{if } d_j(x,y) = g_j(x) - g_j(y) \leq 0, \\
                        	    1, & \mbox{if } d_j(x,y) = g_j(x) - g_j(y)> 0.
                        	\end{array}
                        \right.
                \end{equation}  
                
        It is worth emphasizing that when applying a preferential function regarding a criterion $g_j$ that must be minimized, the order of the calculation of $d_j$ is reversed. Therefore, $P_j (x, y) = P [d_j (x, y)] = P[g_j (y) - g_j (x)].$
    
        In order to apply the Choquet-FlowSort, let us consider set $R_i = R \cup \{a_i\}$, $i=1,\ldots,10$. For a better understanding of the proposed method, we present in detail some computations regarding alternative $a_1$. For the other alternatives, only the results of each step are presented. 
        
        The first step of the FlowSort-Choquet algorithm (S1, in Figure \ref{fig:Algoritmo}) is to conduct the pairwise comparisons between the elements of $R_i$ for each $i=1,\ldots, 10$, applying the chosen preference function, in our case, the preference function of type 1 presented in equation \eqref{eq_pref_function}. Regarding alternative $a_1$, let us consider the set $R_1 = R \cup \{a_1\} = \{r_1, r_2, r_3, r_4, a_1\}$. Equation \eqref{eq_exemplo_a1} presents the computation of the pairwise comparisons between alternative $a_1$ and the other elements of $R_1$ regarding the criterion price, and equation \eqref{eq_exemplo_a12} presents the pairwise comparisons between the other elements of $R_1$ and alternative $a_1$, also regarding the criterion price.
        
        \begin{eqnarray}\label{eq_exemplo_a1}
            d_1(a_1, r_1) = g_1(r_1) - g_1(a_1) = 15,000 - 16,000 \leq 0 &\implies& P_1(a_1, r_1) = 0,\nonumber\\
            d_1(a_1, r_2) = g_1(r_2) - g_1(a_1) = 16,000 - 16,000 \leq 0 &\implies& P_1(a_1, r_2) = 0,\\
            d_1(a_1, r_3) = g_1(r_2) - g_1(a_1) = 17,000 - 16,000 > 0 &\implies& P_1(a_1, r_3) = 1,\nonumber\\
            d_1(a_1, r_3) = g_1(r_2) - g_1(a_1) = 18,000 - 16,000 > 0 &\implies& P_1(a_1, r_4) = 1.\nonumber
        \end{eqnarray}
        \begin{eqnarray}\label{eq_exemplo_a12}
            d_1(r_1,a_1) = g_1(a_1) - g_1(r_1)  = 16,000 - 15,000 > 0 &\implies& P_1(r_1,a_1) = 1,\nonumber\\
            d_1(r_2, a_1) = g_1(a_1) - g_1(r_2) = 16,000 - 16,000 \leq 0 &\implies& P_1(r_2, a_1) = 0,\\
            d_1(r_3, a_1) = g_1(a_1) - g_1(r_2) = 16,000 - 17,000 \leq 0 &\implies& P_1(r_3, a_1) = 0,\nonumber\\
            d_1(r_3, a_1) = g_1(a_1) - g_1(r_2) = 16,000 - 18,000 \leq 0 &\implies& P_1(r_3, a_1) = 0.\nonumber
        \end{eqnarray}
        
        Table \ref{tab:pref_function} presents the value of preference functions of alternative $a_1$ and the other elements of $R_1$ regarding all criteria. Analogously, we can compute the preference functions between each alternative $a_i$ and the limiting profiles of $R_i$ ($P_1(a_i, r_h)~h=1, \ldots, 4, ~i=1, \ldots, 10$), and the preference functions between each limiting profile of $R_i$ and all alternatives $a_i$ ($P_1(r_h, a_i)~h=1, \ldots, 4, ~i=1, \ldots, 10$).
     
        \begin{table}[htb!]
           \centering \small
            \caption{Values of the preference functions between $a_i$ and the limiting profiles of $R_i$ regarding each criteria, i.e, $P_j(a_1, r_h), h=1, \ldots, 4, j=1, \ldots, 4$.}
            \label{tab:pref_function}
            \begin{tabular}{c|cccc|cccc|cccc|cccc}\hline
                 Cars & \multicolumn{4}{c|}{Price} & \multicolumn{4}{c|}{Acceleration} & \multicolumn{4}{c|}{Max Speed} &\multicolumn{4}{c}{Consumption} \\\hline
                	& $r_1$ & $r_2$ & $r_3$	&	$r_4$ & $r_1$ & $r_2$ & $r_3$	&	$r_4$ & $r_1$ & $r_2$ & $r_3$	&	$r_4$ & $r_1$ & $r_2$ & $r_3$	&	$r_4$ \\\hline
                    $a_1$ & 0 & 0 & 1&  1& 0 & 0&1&1&0&1&1&1&0&1&1&1  \\		\hline	
            \end{tabular}
        \end{table}
        
        The second step (S2) is the computation of the Choquet-outranking degrees using equation \eqref{choquet_outeranking_degree2}. We present the computation of $CI_{\pi}(a_1, r_2)$ in equation \eqref{example_Choquet-outranking}. Following the same idea, we obtain $CI_{\pi}(a_1, r_1) = 0$, $CI_{\pi}(a_1, r_3) = 1$ and $CI_{\pi}(a_1, r_4)=1$. Analogously, we can compute the Choquet-outranking degrees of each limiting profile $r_h, h=1, \ldots, 4$ regarding the other elements of $R_i$.
        \begin{eqnarray}
            CI_{\pi}(a_1, r_1) &=& \min\{P_3(a_1,r_2), P_4(a_1,r_2)\} \times I_{34} + \max\{P_2(a_1,r_2),                                               P_3(a_1,r_2)\} \times |I_{23}|\nonumber\\
                               &+& P_1(a_1,r_2) \times I_1 + P_2(a_1,r_2) \times (I_2 - \frac{1}{2}|I_{23}|) \nonumber\\
                               &+& P_3(a_1,r_2) \times (I_3 - \frac{1}{2}\left(|I_{23}|+|I_{34}|)\right) 
                               + P_4(a_1,r_2) \times (I_4 - \frac{1}{2}(|I_{34}|))\\\label{example_Choquet-outranking} 
                              &=& 0.10 + 0.08 + 0.07 + 0.33 = 0.58.  
        \end{eqnarray}
        Now we are able to compute the Choquet-flows of alternatives (S3) and profiles (S4) using equations \eqref{pos_flow_flow-choquet2}, \eqref{neg_flow_flow-choquet2} and \eqref{net_flow_flow-choquet2}. The results of the positive, negative and net Choquet-flows for all alternatives and profiles are presented in Table \ref{tab:flows}. 
        \begin{table}[htb!]
            \centering\small
            \caption{Positive, negative and net Choquet-flows for all alternatives and profiles.}
            \label{tab:flows}
            \begin{tabular}{c|ccccc|ccccc|ccccc}\hline
                    &  \multicolumn{5}{c}{$\phi^{+}_{R_i,CI_{\pi}}$} & \multicolumn{5}{c}{$\phi^{-}_{R_i,CI_{\pi}}$} & \multicolumn{5}{c}{$\phi_{R_i,CI_{\pi}}$} \\\hline
               $i$     &	$a_i$	&	$r_1$	&	$r_2$	&	$r_3$	&	$r_4$	&	$a_i$	&	$r_1$	&	$r_2$	&	$r_3$	&	$r_4$	&	$a_i$	&	$r_1$	&	$r_2$	&	$r_3$	&	$r_4$ \\\hline
            1	&	0.645	&	1	&	0.563	&	0.250	&	0	&	0.313	&	0	&	0.395	&	0.750	&	1	&	0.333	&	1	&	0.168	&	-0.500	&	-1	\\
            2	&	0.375	&	1	&	0.688	&	0.330	&	0	&	0.518	&	0	&	0.313	&	0.563	&	1	&	-0.143	&	1	&	0.375	&	-0.233	&	-1	\\
            3	&	0.518	&	1	&	0.645	&	0.333	&	0	&	0.478	&	0	&	0.375	&	0.643	&	1	&	0.040	&	1	&	0.270	&	-0.310	&	-1	\\
            4	&	0.645	&	1	&	0.580	&	0.250	&	0	&	0.330	&	0	&	0.395	&	0.750	&	1	&	0.315	&	1	&	0.185	&	-0.500	&	-1	\\
            5	&	0.580	&	1	&	0.645	&	0.250	&	0	&	0.395	&	0	&	0.330	&	0.750	&	1	&	0.185	&	1	&	0.315	&	-0.500	&	-1	\\
            6	&	0.563	&	1	&	0.563	&	0.250	&	0	&	0.313	&	0	&	0.313	&	0.750	&	1	&	0.250	&	1	&	0.250	&	-0.500	&	-1	\\
            7	&	0.375	&	1	&	0.708	&	0.250	&	0	&	0.458	&	0	&	0.313	&	0.563	&	1	&	-0.083	&	1	&	0.395	&	-0.313	&	-1	\\
            8	&	0.395	&	1	&	0.750	&	0.375	&	0	&	0.625	&	0	&	0.250	&	0.645	&	1	&	-0.230	&	1	&	0.500	&	-0.270	&	-1	\\
            9	&	0.375	&	1	&	0.708	&	0.313	&	0	&	0.520	&	0	&	0.313	&	0.563	&	1	&	-0.145	&	1	&	0.395	&	-0.250	&	-1	\\
            10	&	0.288	&	1	&	0.750	&	0.458	&	0	&	0.708	&	0	&	0.250	&	0.538	&	1	&	-0.420	&	1	&	0.500	&	-0.080	&	-1	\\ \hline
            \end{tabular}
        \end{table}
         
         As an example of the computation of the Choquet-flow, the calculation of the positive Choquet-flow of alternative $a_1$ is given by
        \begin{eqnarray*}
            \phi^{+}_{R_1,CI_{\pi}}(a_1) &=& 
            \frac{1}{|R_1| - 1}\sum_{y \in \{R_1\} - a_1} CI_{\pi}(a_1, y)\\
                                         &=&\frac{1}{(5 - 1)} \times \left[CI_{\pi}(a_1, r_1)+CI_{\pi}(a_1, r_2)+CI_{\pi}(a_1, r_3)+CI_{\pi}(a_1, r_4)\right]\\
                                         &=&\frac{1}{4} \times(0+0.58+1+1) = 0.645.\\
        \end{eqnarray*}
         The last step (S5) is the assignment of the alternatives to the predefined categories according to the rules \eqref{rule1}, \eqref{rule2} and \eqref{rule1}. The assignments are presented in Table \ref{tab:assignments}.

        \begin{table}[htb!]
            \centering\small
            \caption{Assignments of the ten car models based on the three different assignment rules applying FlowSort-Choquet: assignment rule \eqref{rule1}, based on the positive Choquet-flow, assignment rule 2 \eqref{rule2}, based on the negative Choquet-flow, and assignment rule 3 \eqref{rule3}, based on the net Choquet-flow.}
            \label{tab:assignments}
            \begin{tabular}{c|ccc}\hline
                Cars &  Positive Choquet-flow & Negative Choquet-flow  & Net Choquet-flow\\\hline
                    $a_1$	&	$K_1$	&	$K_1$	&	$K_1$	\\
                    $a_2$	&	$K_2$	&	$K_2$	&	$K_2$	\\
                    $a_3$	&	$K_2$	&	$K_2$	&	$K_2$	\\
                    $a_4$	&	$K_1$	&	$K_1$	&	$K_1$	\\
                    $a_5$	&	$K_2$	&	$K_2$	&	$K_2$	\\
                    $a_6$	&	$K_2$	&	$K_1$	&	$K_2$	\\
                    $a_7$	&	$K_2$	&	$K_2$	&	$K_2$	\\
                    $a_8$	&	$K_2$	&	$K_2$	&	$K_2$	\\
                    $a_9$	&	$K_2$	&	$K_2$	&	$K_2$	\\
                    $a_{10}$	&	$K_3$	&	$K_3$	&	$K_3$	\\\hline
            \end{tabular}
        \end{table}
        
        \section{Numerical experiments and some discussions}
        
        \noindent In this section, we conduct some numerical experiments showing how the consideration or not of either synergy or redundancy between criteria may impact the classification of the alternatives. The numerical experiments are based on the same example presented in Section \ref{sec_example}, but with different capacity values.
        
        Firstly, we apply the FlowSort-Choquet without any interaction information. Then, we introduce either synergy or redundancy between criteria. Moreover, we vary values of the interaction indices in order to identify how the indexes magnitude may also influence the results. To keep the same basis of comparison for all scenarios, the importance of criteria are defined using the interaction and Shaply indexes. The interaction and Shaply indexes assumed in each scenario are presented in Table \ref{tab:scenarios}. For each scenario, we applied the FlowSort-Choquet using equation \eqref{choquet_outeranking_degree2}. The results obtained are presented in Table \ref{tab:assignments_scenarios}.
        
        \begin{table}[htb!]
            \centering \small
            \caption{Scenarios with different values of interaction and Shaply indexes.}
            \label{tab:scenarios}
            \begin{tabular}{c|rrrrrr}\hline
               Scenarios & $I_{1}$ &  $I_{2}$ & $I_{3}$ & $I_{4}$ &  $I_{23}$ & $I_{34}$\\\hline
               Scenario 0 & 0.25 &	0.25 &	0.25 &	0.25 &	- &	-\\\hline
               Scenario 1 & 0.25 &	0.25 &	0.25 &	0.25 &	- &	0.04\\
               Scenario 2 & 0.25 &	0.25 &	0.25 &	0.25 &	- &	0.24\\
               Scenario 3 & 0.25 &	0.25 &	0.25 &	0.25 &	- &	1,0E-15 \\
               Scenario 4 & 0.25 &	0.25 &	0.25 &	0.25 &	- &	9,9E-13 \\
               Scenario 5 & 0.25 &	0.25 &	0.25 &	0.25 &	-0.20 &	-\\
               Scenario 6 & 0.25 &	0.25 &	0.25 &	0.25 &	-0.10 &	0.20\\
               \hline
            \end{tabular}
        \end{table}
        
        As no interaction information is considered in Scenario 0, the application of FlowSort-Choquet is reduced to the application of FlowSort. Therefore, Scenario 0 gives us the opportunity of comparing the proposed method with the traditional FlowSort. In Scenarios 1 e 2, we introduce synergy between criteria 3 and 4. When we consider the synergy index $I_{34}=0.04$ between criteria 3 and 4 (Scenario 1), a change in the results can be sighted: alternative $a_2$ is now assigned to a better category, changing from $K_3$ to $K_2$. This is due the fact that alternative $a_2$ has good evaluations in both maximum speed and consumption criteria. Since a speedy car also having a low consumption is well appreciated, alternative $a_2$ is better classified due to the synergy factor. Increasing the value of this synergy index to $I_{34}=0.24$ (Scenario 2), we obtain the same result as in Scenario 1: alternative $a_2$ changes from $K_3$ to $K_2$. 
        
        In Scenario 3, where a very small synergy value is assumed between criteria 3 and 4 ($I_{34} =$1,0E-15), the result is the same as in Scenario 0, when no synergy value is assumed. However, assuming a synergy index greater than 1,0E-15 (Scenario 4) leads $a_2$ to be assigned to $K_2$, as occurred in Scenarios 2 and 3. Therefore, it can be noticed that different synergy values can lead to the same result, and moreover there is a minimum value that may change the result with any increase in it. Thus, a question arises: which synergy value should be assumed? This is a difficult question to answer but also a non-recent query in MCDM/A \citep{VETSCHERA2017244, Pelissari2019}. Indeed, the elicitation of preference parameters, in our case the elicitation of the capacities, is considered a challenging issue in the application of MCDM/A methods and there are several approaches that can be applied in order to indirectly elicit capacities without the need to define exact values for them. For instance, a stochastic-based approach (\cite{Angilella2015a}) or an incremental preference elicitation approach \citep{BENABBOU2017152}.

        In Scenario 5, we introduce a redundancy between criteria 2 and 3  ($I_{23}=-0.20$). When compared to Scenario 0, alternative $a_6$ is assigned to a better category, changing from $K_2$ to $K_1$. This is due the fact that alternative $a_6$ has relative good evaluations in both acceleration and maximum speed. In spit of alternative $a_3$ has also relative good evaluations in both criteria 2 and 3, it is not assigned to a different category in Scenario 5. This is due to the fact that others factors, such as the limiting profiles adopted, may influence the assignment changing. This means that the Choquet-flow of an alternative can change when considering synergy or redundancy between two criteria, but not necessarily this change is sufficient to cause the alternative to change categories. 

        In Scenario 6, we introduce both synergy between criteria 3 and 4 and redundancy between criteria 2 and 3. The results show a change in the assignment of alternatives $a_3$ and $a_6$, the same ones that experienced changes when considering either synergy or redundancy.
        
        \begin{table}[htb!]
            \centering\small
            \caption{Assignments of the ten car models based on the different scenarios.}
            \label{tab:assignments_scenarios}
            \begin{tabular}{c|c|cccccc}\hline
            Cars    &  Scen0 &  Scen1       &  Scen2       & Scen3  &  Scen4        & Scen5 &  Scen6 \\\hline
            $a_1$	&	$K_1$&	$K_1$	    &	$K_1$      &$K_1$   &$K_1$          &$K_1$          &$K_1$	 \\
            $a_2$	&	$K_3$&  \bm{$K_2$}&\bm{$K_2$}&$K_3$ &\bm{$K_2$}&$K_3$       &\bm{$K_2$}\\
            $a_3$	&	$K_2$&	$K_2$	    &	$K_2$      &$K_2$   &$K_2$          &$K_2$	        &$K_2$	\\
            $a_4$	&	$K_2$&	$K_2$	    &	$K_2$      &$K_2$   &$K_2$          &$K_2$	        &$K_2$	\\
            $a_5$	&	$K_2$&	$K_2$	    &	$K_2$      &$K_2$   &$K_2$          &$K_2$	        &$K_2$	\\
            $a_6$	&	$K_2$&	$K_2$	    &	$K_2$      &$K_2$   &$K_2$          &\bm{$K_1$}     &\bm{$K_1$}	\\
            $a_7$	&	$K_2$&	$K_2$	    &	$K_2$      &$K_2$   &$K_2$          &$K_2$	        &$K_2$	\\
            $a_8$	&	$K_3$&	$K_3$	    &	$K_3$      &$K_3$   &$K_3$          &$K_3$	        &$K_3$	\\
            $a_9$	&	$K_2$&   $K_2$	    &   $K_2$      &$K_2$   &$K_2$          &$K_2$	        &$K_2$\\
            $a_{10}$&	$K_3$&	$K_3$	    &	$K_3$      &$K_3$   &$K_3$          &$K_3$   	    &$K_3$\\\hline
            \end{tabular}
        \end{table}

\section{Conclusion}

    \noindent In this paper, we proposed a new outranking formulation of the Choquet integral for sorting problems in which criteria evaluations can be expressed in heterogeneous scales without any previous conversion into a common scale. Indeed, the common scale required by the Choquet integral is built inside the FlowSort framework by using preference functions, in a way that no prior commensurability assumption on criteria evaluations is required. Moreover, instead of the weighted sum, in FlowSort-Choquet the Choquet integral is applied and, consequently, the proposed method can be seen as an extension of the FlowSort method for interacting criteria. 
    
    We applied the proposed method to an example, showing that it is conceptually simple to be implemented. We also conducted some numerical experiments, which showed us that when synergy or redundancy between criteria are considered, the same alternatives may end up assigned to different categories, when comparing to the situation in which no interaction is taken into account. Moreover, the magnitude of the interaction index also influences the assignments. At the same time, it is difficult for the DM to identify and define interaction values in order to properly represent his/her preferences. Therefore, a natural perspective for future work is the extension of the proposed approach to deal with the elicitation of the capacities.
    
\section*{Acknowledgments}
    \noindent This study was funded by: grant \#2018/23447-4, São Paulo Research Foundation (FAPESP), and by the Brazilian National Council for Scientific and Technological Development (CNPq).
   
\appendix

\renewcommand{\bibname}{Bibliografia} 
\addcontentsline{toc}{chapter}{Bibliografia}
\bibliography{Ref}    

\end{document}